\documentclass[twoside,11pt]{article} %

 \usepackage{hyperref} 
 \usepackage{amsmath}
 \usepackage{amsthm}
 \usepackage{amsfonts}
 \usepackage{amssymb}
 \usepackage{amsopn}
 \usepackage{bbm}
 \usepackage{dsfont}
 \usepackage{mathtools,enumerate,amssymb}
 \usepackage{mathrsfs}
 \usepackage{mathtools}
 \usepackage[T1]{fontenc} 
 \usepackage[utf8]{inputenc} 
 \usepackage[english]{babel} 
 \usepackage{lipsum} 
 \usepackage{url} 

\newcommand{\proofend}{\hfill $\Box$\\}   

\newtheorem{proposition}{\bf Proposition}
\newtheorem{theorem}{\bf Theorem}
\newtheorem{example}{\bf Example}
\newtheorem{remark}{\bf Remark}
\newtheorem{lemma}{\bf Lemma}
\newtheorem{corollary}{\bf Corollary}
\newtheorem{definition}{\bf Definition}

\title{\Large \bf ON THE FRACTIONAL PROBABILISTIC\\[2mm]
 TAYLOR'S AND MEAN VALUE THEOREMS}

 \author{
 {\sc Antonio Di Crescenzo}\footnote{
 Dipartimento di Matematica, Universit\`a degli Studi di Salerno,
 Via Giovanni Paolo II, 132, 84084 Fisciano (SA), Italy, email: adicrescenzo@unisa.it}
 \and 
 {\sc Alessandra Meoli}\footnote{
 Dipartimento di Matematica, Universit\`a degli Studi di Salerno,
 Via Giovanni Paolo II, 132, 84084 Fisciano (SA), Italy, email: ameoli@unisa.it}
}

\date{\normalsize 
\bf First published in {\em Fractional Calculus and Applied Analysis}, \\
Vol.\ 19, n.\ 4, p.\ 921--939 \ \ \copyright\ 2016 by De Gruyter
}

\begin{document}
 
 \maketitle

 \begin{abstract}
In order to develop certain fractional probabilistic analogues of Taylor's theorem and mean value theorem,  
we introduce the $n$th-order fractional equilibrium distribution in terms of the Weyl fractional integral and 
investigate its main properties. Specifically, we show a characterization result by which the $n$th-order 
fractional equilibrium distribution is identical to the starting distribution if and only if it is exponential.
The $n$th-order fractional equilibrium density is then used to prove a fractional probabilistic Taylor's theorem 
based on derivatives of Riemann-Liouville type. A fractional analogue of the probabilistic mean value theorem
is thus developed for pairs of nonnegative random variables ordered according to the survival bounded 
stochastic order. We also provide some related results, both involving the normalized moments and a 
fractional extension of the variance, and a formula of interest to actuarial science. 
In conclusion we discuss the probabilistic Taylor's theorem based on fractional Caputo derivatives. 

 \medskip

{\it MSC 2010\/}: Primary 60E99;
                  Secondary 26A33, 26A24.

 \smallskip

{\it Key Words and Phrases}: characterization of exponential distribution, fractional calculus, fractional equilibrium distribution, 
generalized Taylor's formula, mean value theorem, survival bounded order.

 \end{abstract}

 \section{Introduction and background}\label{sec:1}

Taylor's theorem is the most important result in differential calculus. In fact, given the derivatives of a function at a single point, it provides insight into the behavior of the function at nearby points. Motivated by the large numbers of its applications, researchers have shown a heightened interest in the extensions of this theorem. For instance, Massey and Whitt \cite{MasseyWhitt1993} derived probabilistic generalizations of the fundamental theorem of calculus and Taylor's theorem by making the argument interval random and expressing the remainder terms by means of iterates of the equilibrium residual-lifetime distribution from the theory of stochastic point processes. Lin \cite{Lin1994} modified Massey and Whitt's probabilistic generalization of Taylor's theorem and gave a natural proof by using an explicit form for the density function of high-order equilibrium distribution. In a similar spirit to these probabilistic extensions of Taylor's theorem, Di Crescenzo \cite{DiCrescenzo1999} gave a probabilistic analogue of the mean value theorem. The previous results have direct applications to queueing and reliability theory. However, probabilistic generalizations are not the only ones. Indeed, fractional Taylor series have been introduced with the idea of approximating non-integer power law functions. Here we recall the most interesting ones. Trujillo et al. \cite{Trujillo1999} established a Riemann-Liouville generalized Taylor's formula, in which the coefficients are expressed in terms of the Riemann-Liouville fractional derivative. On the other hand, Odibat et al. \cite{Odibat2007} expressed the coefficients of a generalized Taylor's formula in terms of the Caputo fractional derivatives. 
In the aforementioned papers, an application of the generalized Taylor's formula to the resolution of fractional differential equations is also shown. Great emphasis has been placed on fractional Lagrange and Cauchy type mean value theorems too 
(cf.\ \cite{Guo2012}, \cite{Odibat2007}, \cite{Pecaric2005}, \cite{Trujillo1999}, for example). Inspired by such improvements, in our present investigation we propose to unify these two approaches by presenting a fractional probabilistic Taylor's theorem and a fractional probabilistic mean value theorem.
\par
We start by briefly recalling some basic definitions and properties of fractional integrals and derivatives of Riemann-Liouville and Weyl type as well as some notions on a generalized Taylor's formula that are pertinent to the developments in this paper. For more details on fractional calculus we refer the reader to \cite{Gorenflo2014}.

\subsection{Background on fractional integrals}
Let $ \Omega = \left [ a,b \right ]$  $\left ( -\infty<a<b\leq+\infty \right ) $ be an interval on the real axis $ \mathbb{R} $. The \textit{progressive} or \textit{right-handed fractional integral} $ I_{a^{+}}^{\alpha }f $ of order $ \alpha $ is defined by 
\begin{equation}\label{RL:a+}
I_{a+}^{\alpha }f \left ( x \right ):=\frac{1}{\Gamma \left ( \alpha  \right )}\int_{a}^{x}\left ( x-t \right )^{\alpha -1}f\left ( t \right )dt,\qquad a<x<b,\quad \alpha>0.
\end{equation}
\noindent Here $ \Gamma \left ( \alpha  \right ) $ is the Gamma function and the function $ f(x) $ is assumed to be well-behaved, in order to ensure the finiteness of the values $ I_{a+}^{\alpha }f \left ( x \right ) $ for $ a<x<b $.
\noindent A major property of the fractional integral (\ref{RL:a+}) is the additive index law \textit{(semi-group property)}, according to which
\begin{equation*}
I_{a+}^{\alpha }I_{a+}^{\beta }= I_{a+}^{\alpha+\beta  },\quad 
\alpha ,\beta \geq 0,
\end{equation*}
\noindent where, for complementation, $ I_{a+}^{0}
:=\mathds{I} $ (identity operator).
\par
The \textit{progressive fractional derivative of order} $ \alpha $ is defined by
\begin{equation}\label{RLderivative}
D_{a+}^{\alpha }f(x):=D^{m}I_{a+}^{m-\alpha }f(x),\qquad \text{$ a<x<b,\; m-1<\alpha \leq m $},
\end{equation}
\noindent where $ m $ is a positive integer, and $ D_{a+}^{0}
:=\mathds{I} $ (identity operator). Furthermore, the \textit{sequential fractional derivative} is denoted by
\begin{equation*}
D_{a+}^{n\alpha }=\underbrace{D_{a+}^{\alpha }\dots D_{a+}^{\alpha }}_{n\mathrm{\,times}},
\end{equation*}
\noindent where $ n\in \mathbb{N}\equiv \{0,1,\dots\} $.
\par
A fractional integral over an unbounded interval can also be defined.
Specifically, if the function $ f(x) $ is locally integrable in $ -\infty\leq a<x<+\infty $, and behaves well enough for $x\rightarrow +\infty  $,  the \textit{Weyl fractional integral of order $ \alpha $} is defined as
\begin{equation}\label{Weyl-}
I_{-}^{\alpha }f(x):=\frac{1}{\Gamma \left ( \alpha  \right )}\int_{x}^{+\infty}\left ( t-x \right )^{\alpha -1}f\left ( t \right )dt,\qquad a<x<+\infty,\quad \alpha>0.
\end{equation}
\noindent Also for the Weyl fractional integral the corresponding \textit{semigroup property} holds:
\begin{equation}\label{SemiGroup}
I_{-}^{\alpha }I_{-}^{\beta }= I_{-}^{\alpha +\beta },\quad \alpha ,\beta \geq 0,
\end{equation}
where, again for complementation $ 
I_{-}^{0}:=\mathds{I} $.
\subsection{Background on a generalized Taylor's formula}
Let $ \Omega  $ be a real interval and $ \alpha \in [0,1) $. Let $ F\left( \Omega\right) $ denote the space of Lebesgue measurable
functions with domain in $ \Omega  $ and suppose that $ x_{0}\in \Omega $. Then a function $ f $ is called $ \alpha$-\textit{continuous in} $ x_{0} $ if there exists $ \lambda\in [0,1-\alpha)  $ for which the function $ h $ given by
\begin{equation*}
h\left(x\right)=\left | x-x_{0} \right |^{\lambda }f\left ( x \right )
\end{equation*}
is continuous in $ x_{0} $. Moreover, $ f $ is called $ 1$-\textit{continuous in} $ x_{0} $ if it is continuous in $ x_{0} $, and $ \alpha$-\textit{continuous on }$ \Omega $ if it is $ \alpha$-continuous in $ x $ for every $ x\in\Omega $. We   denote, for convenience, the class of $ \alpha$-continuous functions on $ \Omega $ by $ C_{\alpha}\left(\Omega\right) $, so that $ C_{1}\left(\Omega\right)=C\left(\Omega\right) $.
\par
For $ a\in\Omega $, a function $ f $ is called $ a$-\textit{singular of order $ \alpha $} if
\begin{equation*}
\lim _{x\rightarrow a}\frac{f(x)}{\left | x-a \right |^{\alpha -1}}= k<\infty\quad\mathrm{and}\quad k\neq 0.
\end{equation*}
Let $ \alpha\in\mathbb{R^{+}},\,a\in\Omega$ and let $ E $ be an interval, $ E \subset\Omega $, such that $ a\leq x $ for every $ x\in E$. Then we write
\begin{equation*}
_{a}\mathbf{I}_{\alpha }\left ( E \right )= \left \{ f\in F \left ( \Omega \right ):I_{a+}^{\alpha }f\left ( x \right )\;\textrm{exists and it is finite}\;\forall x\in E\right \}.
\end{equation*}
Recently, Trujillo et al. (cf. Theorem 4.1 of \cite{Trujillo1999}) proved the following result.
\begin{theorem}\label{ThTrujillo}
Set $ \alpha \in \left [ 0,1 \right ] $ and $ n\in \mathbb{N} $. Let $ g $ be a continuous function in $ \left ( a,b  \right ] $ satisfying the following conditions: 
\begin{enumerate}
\item[(i)] $ \forall\,  j= 1,\dots,n, D_{a}^{j\alpha }g\in C\left ( \left ( a,b \right ] \right ) $ and $ D_{a}^{j\alpha }g\in \,_{a}\mathbf{I}_{\alpha }\left ( \left [ a,b \right ] \right ) $;
\item[(ii)] $  D_{a}^{(n+1)\alpha }g $ is continuous on $ \left [ a,b \right ] $;
\item[(iii)] If $ \alpha < 1/2 $ then, for each $ j\in \mathbb{N}, 1\leq j \leq n $, such that $ (j+1)\alpha<1 $, $ D_{a}^{\left ( j+1 \right )\alpha }g\left ( x \right ) $ is $ \gamma$-continuous in $ x=a $ for some $ \gamma $, $ 1-\left ( j+1 \right )\alpha\leq \gamma \leq 1$, or $ a$-singular of order $ \alpha $.
\end{enumerate}
Then, $ \forall x\,\in(a,b] $,
\begin{equation*} 
g(x)=\sum_{j=0}^{n}\frac{c_{j}(x-a)^{\left ( j+1 \right )\alpha -1}}{\Gamma (\left ( j+1 \right )\alpha ) }+R_{n}(x,a),
\end{equation*}
with
\begin{equation*}
R_{n}(x,a)=\frac{D_{a}^{\left ( n+1 \right )\alpha }g\left ( \xi  \right )}{\Gamma \left ( \left ( n+1 \right )\alpha +1 \right )}\left( x-a\right)^{\left ( n+1 \right )\alpha },\qquad a\leq \xi\leq x,
\end{equation*}
and 
\begin{equation*}
c_{j}= \Gamma \left ( \alpha  \right )\left [ \left( x-a\right) ^{1-\alpha } D_{a}^{j\alpha }g(x)\right ]\left ( a^{+} \right )=I_{a}^{1-\alpha }D_{a}^{j\alpha }g\left ( a^{+} \right )
\end{equation*}
for each $ j\in\mathbb{N},\,0\leq j\leq n $.
\end{theorem}
\subsection{Plan of the paper}
The paper is organized as follows. In Section \ref{sec:2}, after recalling the notion of equilibrium distribution, we define a fractional extension of the high-order equilibrium distribution. Then, we give an equivalent version by exploiting the semigroup property of the Weyl fractional integral and derive the explicit expression of the related density function. 
Moreover, by means of the Mellin transform we underline the role played by the fractional equilibrium density in characterizing the exponential distribution. 
In Section \ref{sec:3} we prove a fractional probabilistic Taylor's theorem by using the expression of the $ n $th-order fractional equilibrium density. Section \ref{sec:4} is devoted to the analysis of the fractional analogue of the probabilistic mean value theorem. We first consider pairs of non-negative random variables ordered in a suitable way so as to construct a new random variable, say $ Z_{\alpha} $, which extends the fractional equilibrium operator. 
The fractional probabilistic mean value theorem indeed is given in terms of $ Z_{\alpha} $. 
We also discuss some related results, including a formula of interest to actuarial science.
We stress that all those results are involving the derivatives of Riemann-Liouville type. 
However, in some instances they can be restated also under different setting. 
Indeed, in Section \ref{sec:5} we conclude the paper by exploiting a fractional probabilistic Taylor's theorem 
in the Caputo sense. 
\section{Fractional equilibrium distribution}\label{sec:2}
Let $ X $ be a nonnegative random variable with distribution $ F\left( x\right) = \mathbb{P}\left( X\leq x\right)  $ for $ x\geq 0 $ and with mean $ \mathbb{E}\left[ X\right]\in (0,+\infty)  $. Let $ X_{e} $ be a nonnegative random variable with distribution 
\begin{equation*}
F_{1}\left( x\right) =\mathbb{P}\left ( X_{e}\leq x\right )= \frac{\int_{0}^{x}\overline{F}\left ( y \right )dy}{\mathbb{E}\left[ X\right] },\qquad x\geq 0,
\end{equation*}
or, equivalently, with complementary distribution function
\begin{equation*}
\overline{F}_{1}\left( x\right)=\mathbb{P}\left ( X_{e}> x\right )= \frac{\int_{x}^{+\infty}\overline{F}\left ( y \right )dy}{\mathbb{E}\left[ X\right] },\qquad x\geq 0,
\end{equation*}
where $ \overline{F}=1-F $. The distribution $ F_{1} $ is called \textit{equilibrium distribution with respect to F}. 
Further, suppose $ \mathbb{E}\left[ X^{2}\right]<+\infty  $. Then the equilibrium distribution with respect to $ F_{1} $ is well-defined and it reads
\begin{equation*}
F_{2}\left ( x \right )= \mathbb{P}\left ( X_{e}^{\left ( 2 \right )} \leq x\right )=\frac{\int_{0}^{x}\overline{F}_{1}\left ( y \right )dy}{\mathbb{E}\left [ X_{e} \right ]},\qquad x\geq 0.
\end{equation*}
$ F_{2} $ is known as the second order equilibrium distribution with respect to $ F $. Continuing $ n-2 $ more iterates of this transformation, it is possible to obtain the $ n $th-order equilibrium distribution with respect to $ F $, denoted by $ F_{n} $, provided the required moments of $ X $ are finite.
\par
Hereafter we introduce a fractional version of the $ n $th-order equilibrium distribution. Let $ \alpha\in\mathbb{R}^{+} $ and let $ X $ be a nonnegative random variable with distribution $ F\left( t\right) = \mathbb{P}\left( X\leq t\right)  $ for $ t\geq 0 $ and with moment $ \mathbb{E}\left[ X^{\alpha}\right]\in (0,+\infty)  $. Then we define a random variable $ X_{\alpha}^{(1)} $ whose complementary distribution function is
\begin{equation}\label{FractionalEquilibrium}
\overline{F}_{\,1}^{\,\alpha}\left ( t \right )= \mathbb{P}\left ( X_{\alpha }^{\left ( 1 \right )}>t \right )= \frac{\Gamma \left ( \alpha +1 \right )}{\mathbb{E}\left [ X^{\alpha } \right ]}I_{-}^{\alpha}\overline{F}\left ( t \right ),\qquad t\geq 0,
\end{equation}
where $ I_{-}^{\alpha} $ is the Weyl fractional integral (\ref{Weyl-}) and $ \overline{F}=1-F $. We call the distribution of $ X_{\alpha}^{(1)} $ \textit{fractional equilibrium distribution with respect to $ F $}.
\begin{remark}
Recalling (\ref{Weyl-}), from (\ref{FractionalEquilibrium}) we obtain the following suitable probabilistic interpretation of the distribution function of $ X_{\alpha }^{\left ( 1 \right )} $ in terms of $X$. In fact,
$$
\mathbb{P}\left ( X_{\alpha }^{\left ( 1 \right )}\leq t \right ) 
=\frac{\alpha}{\mathbb{E}\left[ X^{\alpha}\right]}\int_{0}^{+\infty}y^{\alpha -1}\mathbb{P}\left ( y<X\leq y+t \right )dy.
$$
\end{remark}
\par 
Further, suppose $ \mathbb{E}\left[ X^{2\alpha}\right]<+\infty $. Then the \textit{second-order fractional equilibrium distribution with respect to $ F $} is well defined and its complementary distribution function reads 
\begin{equation*}
\overline{F}_{\,2}^{\,\alpha}\left ( t \right )= \mathbb{P}\left ( X_{\alpha }^{\left ( 2 \right )}>t \right )=\frac{\Gamma \left ( 2\alpha +1 \right )}{\Gamma \left ( \alpha +1 \right )}\frac{\mathbb{E}\left [ X^{\alpha } \right ]}{\mathbb{E}\left [ X^{2\alpha } \right ]}I_{-}^{\alpha }\overline{F}_{\,1}^{\,\alpha}\left ( t \right ),\qquad t\geq 0.
\end{equation*}
Generally, we can recursively define the \textit{$ n $th-order fractional complementary equilibrium distribution with respect to $ F $} by
\begin{equation*}
\overline{F}_{\,n}^{\,\alpha }\left ( t \right )=\mathbb{P}\left ( X_{\, \alpha }^{\, \left ( n \right )}>t \right ) = \frac{\Gamma \left ( n\alpha +1 \right )}{\Gamma \left (\left ( n-1 \right )\alpha +1\right )}\frac{\mathbb{E}\left [ X^{\left ( n-1 \right )\alpha } \right ]}{\mathbb{E}\left [X^{n\alpha }  \right ]}\, I_{-}^{\alpha }\overline{F}_{\,n-1}^{\,\alpha }\left( t\right),\qquad t\geq 0,
\end{equation*}
provided that all the moments $ \mathbb{E}[X^{n\alpha}]$, for $ n\in\mathbb{N}$, $n \geq 1$, are finite.
\par
Interestingly enough, $ \overline{F}_{\,n}^{\,\alpha } $ can be alternatively expressed in terms of $ \overline{F} $. Indeed, the following proposition holds.
\begin{proposition}
Let $ \alpha\in\mathbb{R}^{+} $ and let $ X $ be a non-negative random variable with distribution $ F\left( t\right) $ for $ t\geq 0 $. Moreover, suppose that $ \mathbb{E}\left[ X^{n\alpha}\right]\in (0,+\infty)$, with $ n\in\mathbb{N} $, 
$n \geq 1$. Then the \textit{$ n $th-order fractional complementary equilibrium distribution with respect to $ F $} reads
\begin{equation}\label{AlternativeEquilibrium}
\overline{F}_{\,n}^{\,\alpha }\left ( t \right )= \frac{\Gamma \left ( n\alpha +1 \right )}{\mathbb{E}\left [ X^{n\alpha } \right ]}I_{-}^{n\alpha }\overline{F}\left ( t \right ),\qquad t\geq 0.
\end{equation}
\end{proposition}
\proof 
The proof is by induction on $ n $. In fact, when $ n=1 $ formula (\ref{AlternativeEquilibrium}) is true due to definition (\ref{FractionalEquilibrium}). Now let us assume that Eq. (\ref{AlternativeEquilibrium}) holds for some $ n $; then, for $ t\geq 0 $,
\begin{align*}
\overline{F}_{\,n+1}^{\,\alpha }\left ( t \right )&= \frac{\Gamma \left ( (n+1)\alpha +1 \right )}{\Gamma \left (n\alpha +1\right )}\frac{\mathbb{E}\left [ X^{n\alpha } \right ]}{\mathbb{E}\left [X^{(n+1)\alpha }  \right ]}\, I_{-}^{\alpha }\overline{F}_{\,n}^{\,\alpha }\left( t\right)\nonumber \\
&=\frac{\Gamma \left ( (n+1)\alpha +1 \right )}{\Gamma \left (n\alpha +1\right )}\frac{\mathbb{E}\left [ X^{n\alpha } \right ]}{\mathbb{E}\left [X^{(n+1)\alpha }  \right ]}I_{-}^{\alpha}\frac{\Gamma \left ( n\alpha +1 \right )}{\mathbb{E}\left [ X^{n\alpha } \right ]}I_{-}^{n\alpha }\overline{F}\left ( t \right )\nonumber \\
&=\frac{\Gamma \left( \left ( n+1 \right )\alpha +1  \right )}{\mathbb{E}\left [X^{(n+1)\alpha }  \right ]}I_{-}^{\left ( n+1 \right )\alpha }\overline{F}\left ( t \right ).
\end{align*}
\noindent The last equality is valid due to the linearity of the integral and to semigroup property (\ref{SemiGroup}). So the validity of Eq. (\ref{AlternativeEquilibrium}) for $ n $ implies its validity for $ n+1 $. Therefore it is true for all $ n\in\mathbb{N} $, $n\geq 1$.
\proofend 
\par
In order to obtain the explicit expression of the density function of the $ n $th-order fractional equilibrium distribution, we premise the following lemma, which is a generalization of Proposition 4 of \cite{Cheng1999}. 
\begin{lemma}
\label{GeneralizedChengPai}
Let $ X $ be a nonnegative random variable whose moment of order $ n\alpha $ is finite for $ \alpha\in\mathbb{R}^{+} $ and $ n\in\mathbb{N}$, $n\geq 1$. Then
\begin{equation*}
\lim_{x\rightarrow +\infty}\left ( x-t \right )^{n\alpha}\overline{F}\left ( x \right )=0,\quad \forall t\geq 0.
\end{equation*}
\end{lemma}
\proof 
Because the moment of order $ n\alpha $ of $ X $ is finite, we have \begin{equation}\label{FinitenessOfMoment}
\lim_{x\rightarrow +\infty}\int_{x}^{+\infty}y^{n\alpha}dF\left ( y \right )= 0.
\end{equation}
Hence, due to a generalized Markov's inequality,
\begin{equation*}
\lim_{x\rightarrow +\infty}\left ( x-t \right )^{n\alpha}\overline{F}\left ( x \right )\leq \lim_{x\rightarrow +\infty}x^{n\alpha}\overline{F}\left ( x \right )\leq \lim_{x\rightarrow+\infty}\int_{x}^{+\infty}y^{n\alpha}dF\left ( y \right )= 0,
\end{equation*}
this completing the proof.
\proofend
\par
Here and throughout the paper, we denote, for convenience, $ \left ( x \right )^{\alpha-1}_{+}= (x)^{\alpha-1}\mathds{1}_{\{x>0\}} $. The following result concerns the probability density function associated with $ F_{\,n}^{\,\alpha}\left ( t \right ) $.
\begin{proposition}\label{PropositionFractionalDensity}
Let $ X $ be a non-negative random variable with distribution function $ F $ and let $ \mathbb{E}[X^{n\alpha}]<+\infty $ for some integer $ n\geq 1 $ and $ \alpha\in\mathbb{R}^{+} $. Then the density function of $ X_{\, \alpha }^{\, \left ( n \right )} $ is
\begin{equation}\label{FractionalDensity}
f_{n}^{\alpha }\left ( t \right )= \frac{n\alpha \mathbb{E}\left [\left( X-t\right) _{+}^{ n\alpha -1}\right ]}{\mathbb{E}\left [ X^{n\alpha } \right ]},\qquad t\geq 0.
\end{equation}
\end{proposition}
\proof 
By virtue of (\ref{AlternativeEquilibrium}) and (\ref{Weyl-}) we have:
\begin{align*}\label{Stop-loss}
\overline{F}_{\,n}^{\,\alpha }\left ( t \right )&= \frac{\Gamma \left ( n\alpha +1 \right )}{\mathbb{E}\left [ X^{n\alpha } \right ]}I_{-}^{n\alpha }\overline{F}\left ( t \right )\nonumber\\
&=\frac{n\alpha \:  \Gamma \left ( n\alpha  \right )}{\mathbb{E}\left [ X^{n\alpha } \right ]}\frac{1}{\Gamma \left ( n\alpha  \right )}\int_{t}^{+\infty}\left ( x-t \right )^{n\alpha -1}\overline{F}\left ( x \right )dx.\nonumber
\end{align*}
\noindent Due to integration by parts and making use of Lemma \ref{GeneralizedChengPai}, we have:
\begin{align}
\overline{F}_{\,n}^{\,\alpha }\left ( t \right )&=-\frac{1}{\mathbb{E}\left [ X^{n\alpha } \right ]}\int_{t}^{+\infty}\left ( x-t \right )^{n\alpha }d\overline{F}\left ( x \right )\nonumber\\
&=\frac{1}{\mathbb{E}\left [ X^{n\alpha } \right ]}\int_{t}^{+\infty}\left ( x-t \right )^{n\alpha }dF\left ( x \right )\nonumber\\
&=\Bigg(\frac{\mathbb{E}\left [ \left ( X-t \right )_{+}^{n\alpha }\right ]}{\mathbb{E}\left [ X^{n\alpha } \right ]}\Bigg)\\
&= \frac{n\alpha }{\mathbb{E}\left [ X^{n\alpha } \right ]}\int_{t}^{+\infty}dF\left ( x \right )\int_{t}^{x}\left ( x-y \right )^{n\alpha -1}dy\nonumber\\
&=\frac{n\alpha }{\mathbb{E}\left [ X^{n\alpha } \right ]}\int_{t}^{+\infty}dy\int_{y}^{+\infty}\left ( x-y \right )^{n\alpha -1}dF\left ( x \right )\nonumber\\
&=\frac{n\alpha }{\mathbb{E}\left [ X^{n\alpha } \right ]}\int_{t}^{+\infty}\mathbb{E}\left [ \left ( X-y \right )_{+}^{n\alpha -1} \right ]dy,\nonumber
\end{align}
this giving the density function (\ref{FractionalDensity}).
\proofend
\par
We observe that in Proposition \ref{PropositionFractionalDensity} the random variable $ X $ is not necessarily absolutely continuous, unlike $ X_{\, \alpha }^{\, \left ( n \right )} $.
\begin{remark}
Formula (\ref{Stop-loss}) is useful in showing that $ \overline{F}_{\,n}^{\,\alpha }\left ( t \right ) $ is a proper complementary distribution function. Indeed, 
\begin{itemize}
\item[(i)] $ \overline{F}_{\,n}^{\,\alpha }\left ( 0 \right )=\frac{\mathbb{E}\left [ \left ( X-t \right )_{+}^{n\alpha }\right ]}{\mathbb{E}\left [ X^{n\alpha } \right ]}\biggl|_{t=0}=1; $
\item[(ii)] $ \overline{F}_{\,n}^{\,\alpha }\left ( t \right ) $ is decreasing and continuous in $ t\geq 0 $;
\item[(iii)] $ \overline{F}_{\,n}^{\,\alpha }\left ( t \right )\rightarrow 0, $ when $ t\rightarrow +\infty $. In fact, due to (\ref{FinitenessOfMoment}), we have
\begin{equation*}
\lim_{t\rightarrow +\infty}\int_{t}^{+\infty}\left ( x-t \right )^{n\alpha }dF\left ( x \right )\leq \lim_{t\rightarrow +\infty}\int_{t}^{+\infty}x^{n\alpha}dF\left ( x \right )= 0.
\end{equation*} 
\end{itemize} 
\end{remark}
\par 
We now prove a characterization result concerning the fractional equilibrium density (\ref{FractionalDensity}). In fact, if $ X $ is a non-negative random variable with probability density function $ f $, the $ n $th-order fractional equilibrium density associated with $ f $ coincides with $ f $ if and only if $ X $ is exponentially distributed. This extends the well-known result concerning case $ \alpha=1 $.
\begin{theorem}
Let $ X $ be a non-negative random variable with probability density function $ f $. Then, for every $ n\in\mathbb{N} $ and $ \alpha\in\mathbb{R}^{+} $
\begin{equation*}
f_{n}^{\alpha }\left ( t \right )=f(t),\;t\geq 0\qquad\Leftrightarrow\qquad X\sim \mathcal{E}(\lambda),
\end{equation*}
where $ f_{n}^{\alpha }\left ( t \right ) $ is the $ n $th-order fractional equilibrium density (\ref{FractionalDensity}) and $ \mathcal{E}(\lambda) $ is the exponential distribution with parameter $ \lambda\in \mathbb{R}^{+}$.
\end{theorem}
\proof 
First, let us assume that $ X $ is exponentially distributed with parameter $ \lambda $. Since 
\begin{equation*}
\mathbb{E}[X^{n\alpha}]=\frac{\Gamma(n\alpha+1)}{\lambda^{n\alpha}}\quad\mathrm{and}\quad \mathbb{E}\left [ \left( X-t \right )_{+}^{ n\alpha -1}\right]=\lambda^{\-n\alpha}e^{-t\lambda}\Gamma(n\alpha),\quad t\geq 0,
\end{equation*}
by virtue of (\ref{FractionalDensity}) the assertion ``if'' is trivially proved. Conversely, suppose that the density and the $ n $th-order fractional equilibrium density of $ X $ coincide for every $ n\in\mathbb{N} $ and $ \alpha\in\mathbb{R}^{+} $, that is
\begin{equation*}
f_{n}^{\alpha }\left ( t \right )=f(t),\qquad t\geq 0.
\end{equation*}
Due to (\ref{FractionalDensity}), the last equality can be rewritten as
\begin{equation*}
n\alpha \int_{t}^{+\infty}\left ( x-t \right )^{n\alpha -1}f\left ( x \right )dx=f\left ( t \right )\int_{0}^{+\infty}x^{n\alpha }f\left ( x \right )dx,
\end{equation*}
and, on account of (\ref{Weyl-}), as
\begin{equation*}
\Gamma(n\alpha+1)I_{-}^{n\alpha }f(t)=f\left ( t \right )\int_{0}^{+\infty}x^{n\alpha }f\left ( x \right )dx.
\end{equation*} 
Taking the Mellin transform of both sides of this equation yields the functional equation
\begin{equation}\label{TransformedEquation}
\frac{f^{*}\left ( s+n\alpha  \right )}{\Gamma \left (  s+n\alpha \right )}= \frac{f^{*}\left ( n\alpha +1 \right )}{\Gamma \left ( n\alpha +1 \right )}\frac{f^{*}\left ( s \right )}{\Gamma \left ( s \right )},\qquad \Re(s)>0,
\end{equation}
where
\begin{equation*}
f^{*}(s)=\int_{0}^{\infty}x^{s-1}f(x)dx,
\end{equation*}
is the Mellin transform of a function $ f(x) $ (cf. (C.3.21) and (C.3.22) of \cite{Gorenflo2014}). By reducing Eq. (\ref{TransformedEquation}) to a well-known Cauchy equation, we observe that its nontrivial measurable solution (cf. \cite{Hewitt1965} for instance) is
\begin{equation*}
f^{*}(s)=a^{c(s-1)}\Gamma(s)\qquad a>0,\,c\in\mathbb{R}.
\end{equation*} 
By performing the Mellin inversion, we have
\begin{equation*}
f(t)=a^{-c}e^{-a^{-c}t},\qquad t\geq 0.
\end{equation*}
As a consequence, $ X\sim \mathcal{E}(\lambda) $, having set $ \lambda=a^{-c} $, and then the ``only if'' part of the theorem is proved.
\proofend 
\par
In the next proposition, we give the expression of the moments of a random variable following the $ n $th-order fractional equilibrium distribution (\ref{FractionalDensity}). We omit the proof since it is straightforwardly derived after some cumbersome computations.
\begin{proposition}
For $ \alpha\in\mathbb{R}^{+} $ and $ n\in\mathbb{N} $, if $ \mathbb{E}[X^{n\alpha}]<\infty $, then
\begin{equation}\label{moments}
\mathbb{E}\left [\left( X_{\, \alpha }^{\, \left ( n \right )}\right)^{r} \right ]= \frac{n\alpha\,B(n\alpha,r+1) }{\mathbb{E}[X^{n\alpha}]}\, \mathbb{E}[X^{n\alpha+r}],\qquad r\in\mathbb{R}^{+},
\end{equation}
where $ B(x,y) $ is the Beta function.
\end{proposition}
Clearly, when $ \alpha=1 $ the moments (\ref{moments}) identify with the expression for the iterated stationary-excess variables given in the Lemma of Massey and Whitt \cite{MasseyWhitt1993} and in Theorem 2.3 of Harkness and Shantaram \cite{Harkness1969}. 
\section{Fractional probabilistic Taylor's theorem}\label{sec:3}
We now derive a probabilistic generalization of the Riemann-Liouville generalized Taylor's formula shown in Theorem \ref{ThTrujillo}. For convenience, let us denote 
\begin{equation} 
I_{F}=\bigcup_{n=2}^{\infty}\left [ 0,F^{-1}\left ( 1-\frac{1}{n} \right ) \right ],
\label{eq:IF}
\end{equation}
the smallest interval containing both $ {0} $ and the support of the distribution $ F $. Additionally, without loss of generality, we consider the expansion of a function $ g $ about $ t=0 $.
\begin{theorem}\label{ThTaylor}
Let $ 0<\alpha\leq 1 $ and let $ X $ be a nonnegative random variable with distribution $ F $, with moment $ \mathbb{E}[X^{(n+1)\alpha}]<+\infty $ for some integer $ n\geq 0 $ and moments $ \mathbb{E}\left [ X^{\left ( j+1 \right )\alpha -1} \right ]<+\infty  $ for all $ j\in\mathbb{N},0\leq j\leq n $. Suppose that $ g $ is a function defined on $I_{F}$ and 
satisfying the hypoteses (i),(ii) and (iii) of Theorem \ref{ThTrujillo} in $I_{F}$. Assume further $ \mathbb{E}\left [ \left | D_{0}^{\left ( n+1 \right )\alpha }g\left ( X_{\,\alpha }^{\left ( n+1 \right )} \right ) \right |\right ]<+\infty $. Then $ \mathbb{E}\left [ g\left ( X \right ) \right ]<+\infty $ and
\begin{multline}\label{Thesis}
\mathbb{E}\left [ g\left ( X \right ) \right ]= \sum_{j=0}^{n}\frac{c_{j}}{\Gamma \left ( \left ( j+1 \right ) \alpha \right )}\mathbb{E}\left [ X^{\left ( j+1 \right )\alpha -1} \right ] \\
+\frac{\mathbb{E}[X^{(n+1)\alpha}]}{\Gamma \left ( \left ( n+1 \right )\alpha +1 \right )}\,\mathbb{E}\left [ D_{0}^{\left ( n+1 \right )\alpha }g\left ( X_{\,\alpha }^{\left ( n+1 \right )} \right ) \right ],
\end{multline}
with $ c_{j}=\Gamma(\alpha)[x^{1-\alpha}D_{0}^{j\alpha}g(x)](0^{+}) $ for each $ j\in\mathbb{N}, 0\leq j\leq n $, 
where  $ X_{\,\alpha }^{\left ( n+1 \right )} $ has density (\ref{FractionalDensity}).
\end{theorem}
\proof
To begin with, we recall a Riemann-Liouville generalized Taylor's formula with integral remainder term (cf. formula (4.1) of \cite{Trujillo1999}), that is, for $ x\in I_{F} $,
\begin{equation}\label{Taylor}
g(x)=\sum_{j=0}^{n}\frac{c_{j}x^{\left ( j+1 \right )\alpha -1}}{\Gamma (\left ( j+1 \right )\alpha ) }+ R_{n}(x),
\end{equation}
where
\begin{align}\label{resto}
R_{n}(x)&=I_{0}^{\left ( n+1 \right )\alpha }D_{0}^{\left ( n+1 \right )\alpha }g\left ( x \right ) \nonumber\\
&=\frac{1}{\Gamma \left ( \left ( n+1 \right )\alpha  \right )}\int_{0}^{x}\left ( x-t \right )^{\left ( n+1 \right )\alpha -1}D_{0}^{\left ( n+1 \right )\alpha }g\left ( t \right )dt.
\end{align}
Since $ R_{n}(x) $ is continuous on $ I_{F} $, $ R_{n}(X) $ is a true random variable. Therefore, from (\ref{Taylor}) we have
\begin{equation}\label{stella}
\mathbb{E}[g\left ( X \right )]= \sum_{j=0}^{n}\frac{c_{j}}{\Gamma \left ( \left ( j+1 \right ) \alpha \right )}\mathbb{E}[X^{\left ( j+1 \right )\alpha -1}]+\mathbb{E}[R_{n}\left ( X \right )],
\end{equation}
where, from (\ref{resto}),
\begin{equation*}
\mathbb{E}[R_{n}\left ( X \right )]=\frac{1}{\Gamma((n+1)\alpha)}\int_{0}^{+\infty}\int_{0}^{x}D_{0}^{\left ( n+1 \right )\alpha } g\left ( t \right )\left ( x-t \right )^{\left ( n+1 \right )\alpha -1}dt\, dF\left ( x \right ).
\end{equation*}
By making use of Fubini's theorem, the equality above becomes
\begin{equation*}
\mathbb{E}[R_{n}\left ( X \right )]=\frac{1}{\Gamma((n+1)\alpha)}\int_{I_{F}}D_{0}^{\left ( n+1 \right )\alpha }g\left ( t \right )\mathbb{E}\left [ X-t \right ]_{+}^{\left ( n+1 \right )\alpha -1}dt,
\end{equation*}
and in turn, owing to (\ref{FractionalDensity}),
\begin{align}\label{doppiastella}
\mathbb{E}[R_{n}\left ( X \right )]&=\frac{\mathbb{E}[X^{(n+1)\alpha}]}{\left ( n+1 \right )\alpha\, \Gamma \left ( \left ( n+1 \right )\alpha  \right )}\int_{I_{F}}D_{0}^{\left ( n+1 \right )\alpha }g\left ( t \right )f_{n+1}^{\alpha }\left ( t \right )dt\nonumber\\
&=\frac{\mathbb{E}[X^{(n+1)\alpha}]}{ \Gamma \left ( \left ( n+1 \right )\alpha +1 \right )}\mathbb{E}\left [ D_{0}^{\left ( n+1 \right )\alpha }g\left ( X_{\,\alpha }^{\left ( n+1 \right )} \right ) \right ].
\end{align}
Finally, observing that the condition $ \mathbb{E}\left [ \left | D_{0}^{\left ( n+1 \right )\alpha }g\left ( X_{\,\alpha }^{\left ( n+1 \right )} \right ) \right |\right ]<+\infty $ is equivalent to $ \int_{I_{F}} |D_{0}^{\left ( n+1 \right )\alpha }g\left ( t \right ) |\mathbb{E}\left [ X-t \right ]_{+}^{\left ( n+1 \right )\alpha -1}dt<+\infty $, and making use of (\ref{stella}) and (\ref{doppiastella}), the proof of (\ref{Thesis}) is thus completed.
\proofend 
\par
Equation (\ref{Thesis}) can be seen as a fractional version of the probabilistic generalization of Taylor's theorem studied in \cite{Lin1994} and in \cite{MasseyWhitt1993}.
\begin{remark}\label{Oss}
We observe that for $ n=0 $ formula (\ref{Thesis}) becomes 
\begin{equation}\label{RemarkTaylor}
\mathbb{E}\left [ g\left ( X \right ) \right ]=\frac{c_{0}}{\Gamma \left ( \alpha  \right )}\mathbb{E}\left [ X^{\alpha -1} \right ]+\frac{\mathbb{E}\left [ X^{\alpha} \right ]}{\Gamma \left ( \alpha  +1\right )}\mathbb{E}\left [ D_{0}^{\alpha }g\left ( X_{\alpha }^{\left ( 1 \right )} \right ) \right ],
\end{equation}
with $ c_{0}=\Gamma \left ( \alpha  \right )\left [ x^{1-\alpha} g(x)\right ](0^{+}) $, this 
being useful to prove Theorem \ref{Lagrange} below. 
\end{remark}
In the recent years much attention has been paid to the study of the fractional moments of distributions. 
See, for instance, \cite{Mikoschetal2013} and references therein. Motivated by this, in the next corollary we consider the case 
$ g(x)=x^{\beta},\,\beta\in\mathbb{R} $. From Theorem \ref{ThTaylor} we have the following result.
\begin{corollary}
Let $ 0<\alpha\leq 1 $, $ \beta\geq\alpha $ and $ n\leq \frac{\beta-\alpha}{\alpha},\, n\in \mathbb{N} $. Moreover, let $ X $ be a nonnegative random variable with distribution $ F $, with moment $ \mathbb{E}[X^{(n+1)\alpha}]<+\infty $ and moments $ \mathbb{E}[X^{(j+1)\alpha-1}]<+\infty $ for all $ j\in\mathbb{N}, 0\leq j\leq n $. Assume further $ \mathbb{E}\left [ \left | D_{0}^{\left ( n+1 \right )\alpha }\left ( X_{\,\alpha }^{\left ( n+1 \right )} \right )^{\beta} \right |\right ]<+\infty $. Then  $ \mathbb{E}[X^{\beta}]<+\infty $ and
\begin{equation*}
\mathbb{E}[X^{\beta}]=\frac{\mathbb{E}[X^{(n+1)\alpha}]}{\Gamma \left ( \left ( n+1 \right )\alpha +1 \right )}\,\frac{\Gamma(1+\beta)}{\Gamma(1-(n+1)\alpha + \beta)}\mathbb{E}\left [ (X_{\,\alpha }^{\left ( n+1 \right )})^{\beta-(n+1)\alpha}\right],
\end{equation*}
where  $ X_{\,\alpha }^{\left ( n+1 \right )} $ has density (\ref{FractionalDensity}).
\end{corollary}
\proof
Since, in general, for $ k\in\mathbb{N} $
\begin{equation*}
D_{0}^{ k \alpha }x^{\beta }=\frac{\Gamma(1+\beta)}{\Gamma(1-k\alpha+\beta)}x^{\beta-k\alpha},
\end{equation*}
we have
\begin{equation*}
c_{j}=\Gamma(\alpha)\left[\frac{\Gamma(1+\beta)}{\Gamma(1-j\alpha+\beta)}x^{1-(j+1)\alpha+\beta}\right](0^{+})=0, \qquad 0\leq j\leq n.
\end{equation*}
Furthermore, assumption $ \mathbb{E}\left [ \left | D_{0}^{\left ( n+1 \right )\alpha }\left ( X_{\,\alpha }^{\left ( n+1 \right )} \right )^{\beta} \right |\right ]<+\infty $ ensures the finiteness of $ \mathbb{E}\left [ (X_{\,\alpha }^{\left ( n+1 \right )})^{\beta-(n+1)\alpha}\right] $.
Therefore, formula (\ref{Thesis}) reduces to the sole remainder term, and hence the thesis.
\proofend
\section{Fractional probabilistic mean value theorem}\label{sec:4}
In this section we develop the probabilistic analogue of a fractional mean value theorem. 
To this purpose we first recall some stochastic orders and 
introduce a relevant random variable, $ Z_{\alpha} $. 
\par
Let $ X $ be a random variable with distribution function $ F_{X} $ and let $ a=\inf \left \{ x|F_{X}\left ( x \right )>0 \right \} $ and $ b=\sup \left \{ x|F_{X}\left ( x \right )<1 \right \} $. We set for every real $ \alpha>0 $ 
\begin{equation*}
F_{X}^{\left (\alpha  \right ) }\left ( t \right )=
\left\{
\begin{array}{ll}
\dfrac{\mathbb{E}\left [ \left(t-X \right)_{+}^{\alpha -1}\right ]}{\Gamma \left ( \alpha  \right )} &\mbox{if } t>a\\
0& \text{if }t\leq a,
\end{array}
\right.
\end{equation*}
and 
\begin{equation*}
\overline{F}_{X}^{\left (\alpha  \right ) }\left ( t \right )=
\left\{
\begin{array}{ll}
\dfrac{\mathbb{E}\left [\left( X-t\right)_{+}^{\alpha -1} \right ]}{\Gamma \left ( \alpha  \right )}&\mbox{if } t<b\\
0& \text{if }t\geq b.
\end{array}
\right.
\end{equation*}
\par
Ortobelli et al. \cite{Ortobelli2008} define a stochastic order as follows:
\begin{definition}
Let $ X $ and $ Y $ be two random variables. For every $ \alpha>0 $, $ X $ dominates $ Y $ with respect to the $ \alpha $-\textit{bounded stochastic dominance order} $ (X\overset{b}{\underset{\alpha}{\geq}}Y) $ if and only if $ F_{X}^{\left (\alpha  \right ) }\left ( t \right )\leq F_{Y}^{\left (\alpha  \right ) }\left ( t \right ) $ for every $ t $ belonging to $ \mathit{supp}\{X,Y\}\equiv \left [ a,b \right ] $, where $ a,b\in\mathbb{\overline{R}} $ and $ a=\inf\left \{ x|F_{X} (x)+F_{Y}(x)>0\right \} $, $ b=\sup\left \{ x|F_{X} (x)+F_{Y}(x)<2\right \} $.
\end{definition}
Similarly, Ortobelli et al. \cite{Ortobelli2008} define a \textit{survival bounded order} as follows:
\begin{definition}\label{SurvivalBoundedOrder}
For every $ \alpha>0 $, we write $ X \overset{a}{\underset{\textit{sur $ \alpha $}}{\geq}} Y$ if and only if $\overline{F}_{X}^{\left (\alpha  \right ) }\left ( t \right )\leq\overline{F}_{Y}^{\left (\alpha  \right ) }\left ( t \right )  $ for every $ t $ belonging to $ \mathit{supp}\{X,Y\}$.
\end{definition}
We remark that certain random variables cannot be compared with respect to these orders. 
For example, Ortobelli et al. \cite{Ortobelli2008} proved that for any pair of bounded (from above or/and from below) random variables $ X $ and $ Y $ that are continuous on the extremes of their support, there is no $ \alpha\in (0,1) $ such that $ F_{X}^{\left (\alpha  \right ) }\left ( t \right )\leq F_{Y}^{\left (\alpha  \right ) }\left ( t \right ) $ for all $ t\in\mathit{supp}\{X,Y\} $. However, although $ \alpha $-bounded orders with $ \alpha\in (0,1) $ are not applicable in many cases, they could be useful to rank truncated variables and financial losses, thus resulting of interest from a financial point of view.
\par
We outline that for $ \alpha>1 $ the survival bounded order given in Definition \ref{SurvivalBoundedOrder} is equivalent to the extension to all real $ \alpha>0 $ of the order $ \leq_{c}^{\alpha} $ defined in 1.7.1 of \cite{Stoyan1983}. Moreover, when $ \alpha=2 $, it is equivalent to the increasing convex order $ \leq_{icx} $ (cf. Section 4.A of Shaked and Shantikumar \cite{Shaked2007}).
\par
The following result comes straightforwardly.
\begin{proposition}\label{ZetaAlfa}
Let $ X $ and $ Y $ be non-negative random variables such that $ \mathbb{E}\left[ X^{\alpha}\right]<\mathbb{E}\left[ Y^{\alpha}\right]<+\infty $ for some $ \alpha>0 $. Then
\begin{equation}\label{Transformation}
f_{Z_{\alpha }}\left ( t \right )=\alpha \frac{\mathbb{E}\left [ \left ( Y-t \right )_{+}^{\alpha -1} \right ]-\mathbb{E}\left [ \left ( X-t \right )_{+}^{\alpha -1} \right ]}{\mathbb{E}\left [ Y^{\alpha } \right ]-\mathbb{E}\left [ X^{\alpha } \right ]},\qquad t\geq 0,
\end{equation}
is the probability density of an absolutely continuous non-negative random variable, say $ Z_{\alpha} $, if and only if $ X \overset{0}{\underset{\textit{sur $ \alpha $}}{\geq}} Y$.
\end{proposition}
We remark that condition $ X \overset{0}{\underset{\textit{sur $ \alpha $}}{\geq}} Y $ ensures that $ \mathbb{E}[X^{\alpha}]\leq \mathbb{E}[Y^{\alpha}] $ for $ \alpha>0 $. Moreover, it is interesting to note that $ Z_{\alpha} $ is necessarily absolutely continuous, in contrast with $ X $ and $ Y $.
\begin{example}
Let $ X $ and $ Y $ be exponential random variables having means $ \mu _{X} $ and $ \mu_{Y} $, $ \mu_{Y}>\mu _{X} >0$, and let $ \alpha\geq 1 $. From (\ref{Transformation}) we obtain the following expression for
the density of $ Z_{\alpha } $:
\begin{equation*}
f_{Z_{\alpha }}\left ( t \right )= \frac{\mu _{Y}^{\alpha -1}e^{-\frac{t}{\mu _{Y}}}-\mu _{X}^{\alpha -1}e^{-\frac{t}{\mu _{X}}}}{\mu _{Y}^{\alpha }-\mu _{X}^{\alpha }},\qquad t\geq 0.
\end{equation*}
\end{example}
\begin{example}
Let $ Y $ be a random variable taking values in $ [0,b) $, with $ b\in(0,+\infty] $, and let $ \mathbb{E}[Y^{\alpha}] $ and $ \mathbb{E}[(Y-t)_{+}^{\alpha-1}] $ be finite, $ 0\leq t\leq b $. Furthermore, we define a random variable $ X $ with cumulative distribution function  
\begin{equation*}
F_{X}(x):=\left\{
\begin{array}{ll}
0, &  x<0, \\
p+(1-p)F_{Y}(x), &  0\leq x\leq b,  \\
1, &  x\geq b,
\end{array}
\right.
\end{equation*}
where $ F_{Y}(x) $ is the cumulative distribution function  of $ Y $ and $ 0<p\leq 1 $. We remark that $ X $ can be viewed as a $ 0 $-inflated version of $ Y $, i.e. $ X=I\cdot Y $, where $ I $ is a Bernoulli r.v. independent of Y. It is easily ascertained that 
\begin{equation*}
f_{Z_{\alpha }}\left ( t \right )=\frac{\alpha\mathbb{E}\left [ \left ( Y-t \right )_{+}^{\alpha -1} \right ] }{\mathbb{E}\left [ Y^{\alpha } \right ]}\equiv f_{Y_{\alpha }^{\left ( 1 \right )}}(t)\equiv f_{X_{\alpha }^{\left ( 1 \right )}}(t),\quad t\geq 0.
\end{equation*}
\end{example}
We note that the density of $ Z_{\alpha } $ given in (\ref{Transformation}) is related to the densities of the fractional equilibrium variables $  X_{\alpha}^{(1)} $ and $  Y_{\alpha}^{(1)} $ which, by virtue of (\ref{FractionalDensity}), are respectively given by
\begin{equation*}
f_{X_{\alpha }^{\left ( 1 \right )}}(t)= \frac{\alpha\mathbb{E}\left [ \left ( X-t \right )_{+}^{\alpha -1} \right ] }{\mathbb{E}\left [ X^{\alpha } \right ]}\quad\mathrm{and}\quad f_{Y_{\alpha }^{\left ( 1 \right )}}(t)= \frac{\alpha\mathbb{E}\left [ \left ( Y-t \right )_{+}^{\alpha -1} \right ] }{\mathbb{E}\left [ Y^{\alpha } \right ]},\quad t\geq 0.
\end{equation*}
\par
Indeed, from (\ref{Transformation}) the following generalized mixture holds:
\begin{equation}\label{Non-convexCombination}
f_{Z_{\alpha }}\left ( t \right )= cf_{Y_{\alpha }^{\left ( 1 \right )}}\left ( t \right )+\left ( 1-c \right )f_{X_{\alpha }^{\left ( 1 \right )}}\left ( t \right ),
\end{equation}
where 
\begin{equation}\label{Coefficient}
c=\frac{\mathbb{E}\left [ Y^{\alpha } \right ]}{\mathbb{E}\left [ Y^{\alpha } \right ]-\mathbb{E}\left [ X^{\alpha } \right ]}\geq 1.
\end{equation}
Such representation is useful to find an expression for the moments of $ Z_{\alpha} $. In fact, from (\ref{moments}), (\ref{Non-convexCombination}) and (\ref{Coefficient}), Proposition \ref{MomentsZetaAlpha} follows immediately.
\begin{proposition}\label{MomentsZetaAlpha}
Let $ \alpha \in \mathbb{R}^{+}$ and suppose that X and Y are two non-negative random variables such that $ \mathbb{E}\left[ X^{\alpha}\right]<\mathbb{E}\left[ Y^{\alpha}\right]<+\infty $, and $ X \overset{0}{\underset{\textit{sur $ \alpha $}}{\geq}} Y $. Then 
\begin{equation}\label{momentsZalpha}
\mathbb{E}\left [ Z_{\alpha } ^{r}\right ]= \frac{\alpha B\left ( \alpha ,r+1 \right )}{\mathbb{E}\left [ Y^{\alpha } \right ]-\mathbb{E}[X^{\alpha }]}\left \{ \mathbb{E}\left [ Y^{\alpha+r } \right ]-\mathbb{E}[X^{\alpha+r }] \right \},\qquad r\in\mathbb{R}^{+}.
\end{equation}
\end{proposition}
\par 
With the notation of 1.C(3) of \cite{Marshall2007}, let $ \lambda_{\alpha}(X) $ denote the normalized moment of a random variable $ X $, i.e.
\begin{equation}\label{NormalizedMoments}
\lambda_{\alpha} (X)=\frac{\mathbb{E}[X^{\alpha}]}{\Gamma(\alpha+1)},\qquad \alpha>0.
\end{equation}
\par 
We are now ready to prove the main result of this section.
\begin{theorem}\label{Lagrange}
Let $ 0<\alpha\leq 1 $. Suppose that X and Y are two non-negative random variables such that $ \mathbb{E}\left[ X^{\alpha}\right]<\mathbb{E}\left[ Y^{\alpha}\right]<+\infty $, and $ X \overset{0}{\underset{\textit{sur $ \alpha $}}{\geq}} Y $. Moreover, let Theorem \ref{ThTaylor} hold for some function $ g $. Then
\begin{multline}\label{FracLagrange}
\mathbb{E}\left [ g\left ( Y \right ) \right ]-\mathbb{E}\left [ g\left ( X \right ) \right ]=  \frac{c_{0}}{\Gamma \left ( \alpha  \right )}\left \{ \mathbb{E}\left [ Y^{\alpha -1} \right ] -\mathbb{E}\left [ X^{\alpha -1} \right ] \right \}\\
+\{\lambda_{\alpha}(Y)-\lambda_{\alpha}(X)\}\mathbb{E}\left [ D_{0}^{\alpha }g\left ( Z_{\alpha } \right ) \right ],
\end{multline}
where $ Z_{\alpha } $ is a random variable whose density is defined in (\ref{Transformation}), and $ c_{0}=\Gamma \left ( \alpha  \right )\left [ x^{1-\alpha} g(x)\right ](0^{+}) $.
\end{theorem}
\proof
By applying Theorem \ref{ThTaylor} for $ n=0 $, cf. formula (\ref{RemarkTaylor}), we have
\begin{multline}
\mathbb{E}\left [ g\left ( Y \right ) \right ]-\mathbb{E}\left [ g\left ( X \right ) \right ]= \frac{c_{0}}{\Gamma \left ( \alpha  \right )}\left \{ \mathbb{E}\left [ Y^{\alpha -1} \right ] -\mathbb{E}\left [ X^{\alpha -1} \right ] \right \}\\
+\frac{1}{\Gamma \left ( \alpha +1 \right )}\left \{ \mathbb{E}\left [ Y^{\alpha } \right ]\mathbb{E}\left [ D_{0}^{\alpha } g\left ( Y_{\alpha }^{\left ( 1 \right )} \right )\right ]-\mathbb{E}\left [ X^{\alpha } \right ]\mathbb{E}\left [ D_{0}^{\alpha } g\left ( X_{\alpha }^{\left ( 1 \right )} \right )\right ] \right \}.\nonumber
\end{multline}
From (\ref{Non-convexCombination}), (\ref{Coefficient}) and (\ref{NormalizedMoments}) the theorem is straightforwardly proved.
\proofend
\par
Under certain hypotheses, the Lagrange's Theorem guarantees the existence of a mean value belonging to the interval of interest. With regard to Theorem \ref{Lagrange}, one might therefore expect that a probabilistic analogue of this relation holds too. However, the relation $ X^{\alpha}\leq_{st}Z_{\alpha}\leq_{st} Y^{\alpha} $ does not hold in general. It can be satisfied only when $ \mathbb{E}[X^{\alpha}]\leq\mathbb{E}[Z_{\alpha}]\leq\mathbb{E}[Y^{\alpha}] $, which is case $ (ii) $ of the next Proposition. For simplicity's sake, if $ X $ is a random variable with $ \mathbb{E}[X^{\alpha+1}]<+\infty $, we set
\begin{equation}\label{FractionalVariance}
V_{\alpha }\left ( X \right ):= \mathbb{E}\left [ X^{\alpha +1} \right ]-\alpha \left ( \mathbb{E} \left [ X^{\alpha} \right ]\right )^{2}, \qquad \alpha\in\mathbb{R}^{+},
\end{equation}
which turns out to be a fractional extension of the variance of $ X $.
\begin{proposition}
Let $ 0<\alpha\leq 1 $ and let $ X $ and $ Y $ satisfy the assumptions of Theorem \ref{Lagrange}, with $ \mathbb{E}[X^{\alpha+1}] $ and $ \mathbb{E}[Y^{\alpha+1}] $ finite. Then,\smallskip
\begin{enumerate} 
\item[(i)] $ \begin{multlined}[t]
           \mathbb{E}[Z_{\alpha}]\leq\mathbb{E}[X^{\alpha}]\\ \Leftrightarrow\;V_{\alpha }\left ( Y \right )-V_{\alpha }\left ( X \right )\leq -\left \{\mathbb{E}[Y^{\alpha }] -\mathbb{E}[X^{\alpha}]\right \}\left \{ \alpha \mathbb{E}[Y^{\alpha }] -\mathbb{E}[X^{\alpha}] \right \};
      \end{multlined}$
\smallskip     
\item[(ii)]$ \begin{multlined}[t]
            \mathbb{E}[X^{\alpha}]\leq\mathbb{E}[Z_{\alpha}]\leq\mathbb{E}[Y^{\alpha}]\\ \Leftrightarrow\;\left\{\begin{array}{lr}
            V_{\alpha }\left ( Y \right )-V_{\alpha }\left ( X \right )\geq -\left \{\mathbb{E}[Y^{\alpha }] -\mathbb{E}[X^{\alpha}]\right \}\left \{ \alpha \mathbb{E}[Y^{\alpha }] -\mathbb{E}[X^{\alpha}] \right \} \\ 
           V_{\alpha }\left ( Y \right )-V_{\alpha }\left ( X \right )\leq \left \{\mathbb{E}[Y^{\alpha }] -\mathbb{E}[X^{\alpha}]\right \}\left \{\mathbb{E}[Y^{\alpha }] -\alpha \mathbb{E}[X^{\alpha}]\right \};
            \end{array}\right.
            \end{multlined}$\smallskip
\item[(iii)]$\begin{multlined}[t]
                  \mathbb{E}[Y^{\alpha}]\leq \mathbb{E}[Z_{\alpha}]\\ \Leftrightarrow\;V_{\alpha }\left ( Y \right )-V_{\alpha }\left ( X \right )\geq \left \{\mathbb{E}[Y^{\alpha }] -\mathbb{E}[X^{\alpha}]\right \}\left \{\mathbb{E}[Y^{\alpha }] -\alpha \mathbb{E}[X^{\alpha}]\right \};
                  \end{multlined}$
                 \smallskip
                  \item[(iv)]$\begin{multlined}[t] \mathbb{E}[Z_{\alpha }]= \dfrac{2\alpha }{\alpha +1}\, \dfrac{\mathbb{E}[X^{\alpha }]+\mathbb{E}[Y^{\alpha }]}{2}\;\Leftrightarrow \; V_{\alpha }(Y)=V_{\alpha }(X)\end{multlined} $,
            \end{enumerate}\smallskip
where $ V_{\alpha} $ has been defined in (\ref{FractionalVariance}).
\end{proposition}
\proof
It easily follows from the identity 
\begin{equation*}
\frac{\mathbb{E}\left [ Z_{\alpha } \right ]-\mathbb{E}\left [ X^{\alpha } \right ]}{\mathbb{E}[Y^{\alpha }]-\mathbb{E}\left [ X^{\alpha } \right ]}= \frac{1}{\alpha +1}\left \{ \frac{\alpha \mathbb{E}\left [ Y^{\alpha } \right ]-\mathbb{E}\left [ X^{\alpha } \right ]}{\mathbb{E}\left [ Y^{\alpha } \right ]-\mathbb{E}\left [ X^{\alpha } \right ]} +\frac{V_{\alpha }(Y)-V_{\alpha }(X)}{\left ( \mathbb{E}\left [ Y^{\alpha } \right ]-\mathbb{E}\left [ X^{\alpha } \right ] \right )^{2}}\right \},
\end{equation*}
which is a consequence of (\ref{momentsZalpha}) written for $ r=1 $.
\proofend
\begin{corollary}\label{PreApplication}
Let $ 0<\alpha\leq 1 $. Suppose that X and Y are two non-negative random variables such that $ \mathbb{E}\left[ X^{\alpha}\right]<\mathbb{E}\left[ Y^{\alpha}\right]<+\infty $, and $ X \overset{0}{\underset{\textit{sur $ \alpha $}}{\geq}} Y $. Moreover, if $ \beta>\alpha-1 $, let Theorem \ref{ThTaylor} hold for some function $ g(x)\sim x^{\beta}$, $ x\rightarrow 0^{+}$. Then
\begin{equation*}
\mathbb{E}\left [ g\left ( Y \right ) \right ]-\mathbb{E}\left [ g\left ( X \right ) \right ]=\left \{ \lambda_{\alpha}(Y)-\lambda_{\alpha}(X) \right \}\mathbb{E}\left [ D_{0}^{\alpha }g\left ( Z_{\alpha } \right ) \right ],
\end{equation*}
where $ Z_{\alpha } $ is a random variable whose density is defined in (\ref{Transformation}).
\end{corollary}
\proof
We observe that the first term in formula (\ref{FracLagrange}) vanishes, since $ c_{0}\sim x^{1-\alpha+\beta}\big |_{x=0^{+}}=0$, and hence the thesis.
\proofend
\par
As application, we now show a result of interest to actuarial science. 
A \textit{deductible} is a treshold amount, denoted $ d $, which must be exceeded by a loss in order for a claim to be paid. If $ X $ is the severity random variable representing a single loss event, $ X>0 $, and if the deductible is exceeded (that is, if $ X>d $), then the amount paid is $ X-d $. Therefore, for $ d>0 $, the claim amount random variable $ X_{d} $ is defined to be
\begin{equation}\label{PaymentPerLoss}
X_{d}:=(X-d)_{+}=\left\{
\begin{array}{ll}
0,&\text{for $ X\leq d $,}\\
X-d,&\text{for $ X>d $}.
\end{array}
\right.
\end{equation}   
It is clear from equation (\ref{PaymentPerLoss}) that $ X_{d} $ has a mixed distribution. In particular, such random variable has an atom at zero representing the absence of payment because the loss did not exceed $ d $. The interested reader is referred to \cite{Kellison2011} and \cite{Klugman2012} for further information. Let $ b=\sup \left \{ x|F_{X}\left ( x \right )<1 \right \} $. Bearing in mind Definition \ref{ZetaAlfa}, the next Proposition immediately follows from Corollary \ref{PreApplication}.
\begin{proposition}\label{Appication}
Let $ 0<\alpha\leq 1 $ and $ 0<r<s<b $. With reference to (\ref{PaymentPerLoss}), suppose that $ X_{s} \overset{0}{\underset{\textit{sur $ \alpha $}}{\geq}} X_{r} $ and $ \lambda_{\alpha} (X_{s})<\lambda_{\alpha} (X_{r})<+\infty$. Moreover, let $ g $ satisfy the assumptions of Corollary \ref{PreApplication}. Then
\begin{equation*}
\mathbb{E}\left [ g\left ( X_{r} \right ) \right ]-\mathbb{E}\left [ g\left ( X_{s} \right ) \right ]=[\lambda_{\alpha}(X_{r})-\lambda_{\alpha}(X_{s})]\mathbb{E}\left [ D_{0}^{\alpha }g\left ( Z_{\alpha } \right ) \right ],
\end{equation*}
where $ Z_{\alpha } $ is a random variable with density
\begin{equation*}
f_{Z_{\alpha}}(z)=\alpha \frac{\mathbb{E}\left [ \left ( X_{r}-z \right )_{+}^{\alpha -1} \right ]-\mathbb{E}\left [ \left ( X_{s}-z \right )_{+}^{\alpha -1} \right ]}{\mathbb{E}\left [ X_{r}^{\alpha } \right ]-\mathbb{E}\left [ X_{s}^{\alpha } \right ]},\qquad z\geq 0.
\end{equation*}
\end{proposition}
We conclude this section with the following example.
\begin{example}
Let $ 0<\alpha\leq 1 $ and $ 0<r<s $.
\begin{enumerate}
\item[(i)] Let $ X $ be an exponential random variable with parameter $ \lambda $. Owing to (\ref{NormalizedMoments}) and Proposition \ref{Appication}, we have
\begin{equation*}
\mathbb{E}\left [ g\left ( X_{r} \right ) \right ]-\mathbb{E}\left [ g\left ( X_{s} \right ) \right ]=\lambda^{-\alpha}\left(e^{-\lambda r}-e^{-\lambda s}\right)\mathbb{E}\left [ D_{0}^{\alpha }g\left ( Z\right ) \right ],
\end{equation*}
where $ Z $ turns out to be exponentially distributed with parameter $ \lambda $ as well. In this case, it is interesting to note that for $ 0<r<s $ and $ 0<u<v $ it results:
\begin{equation*}
\frac{\mathbb{E}\left [ g\left ( X_{r} \right ) \right ]-\mathbb{E}\left [ g\left ( X_{s} \right ) \right ]}{\mathbb{E}\left [ g\left ( X_{u} \right ) \right ]-\mathbb{E}\left [ g\left ( X_{v} \right ) \right ]}=\frac{e^{-\lambda r}-e^{-\lambda s}}{e^{-\lambda u}-e^{-\lambda v}},
\end{equation*}
which is independent of $ g $.
\item[(ii)] Now let $ X $ be a 2-phase hyperexponential random variable with phase probabilities $ p $ and $ 1-p $, $ 0<p<1 $, and rates $ \lambda_{1} $ and $ \lambda_{2} $. Similarly, it holds
\begin{multline*}
\mathbb{E}\left [ g\left ( X_{r} \right ) \right ]-\mathbb{E}\left [ g\left ( X_{s} \right ) \right ]=\left \{  p\lambda_{1}^{-\alpha}\left(e^{-\lambda_{1} r}-e^{-\lambda_{1} s}\right)\right.\\
\left. +\left(1-p\right)\lambda_{2}^{-\alpha}\left(e^{-\lambda_{2} r}-e^{-\lambda_{2} s}\right)\right \} \mathbb{E}\left [ D_{0}^{\alpha }g\left ( Z_{\alpha } \right ) \right ],
\end{multline*}
where, the density of $ Z_{\alpha } $ is, for $ z\geq 0 $, 
\end{enumerate}
\begin{equation*}
f_{Z_{\alpha }}\left ( z \right )=\frac{p\lambda_{1}^{1-\alpha}e^{-\lambda_{1}z}\left(e^{-\lambda_{1} r}-e^{-\lambda_{1} s}\right)+\left(1-p\right)\lambda_{2}^{1-\alpha}e^{-\lambda_{2}z}\left(e^{-\lambda_{2} r}-e^{-\lambda_{2} s}\right)}{p\lambda_{1}^{-\alpha}\left(e^{-\lambda_{1} r}-e^{-\lambda_{1} s}\right)+\left(1-p\right)\lambda_{2}^{-\alpha}\left(e^{-\lambda_{2} r}-e^{-\lambda_{2} s}\right)}.
\end{equation*}
\end{example}
\section{Concluding remarks}\label{sec:5}
The overall aim of this research is to present a novel Taylor's theorem from a probabilistic and a fractional perspective at the same time and to discuss other related findings. It is meaningful to note that, while the coefficients of our formula (\ref{Thesis}) are expressed in terms of the Riemann-Liouville fractional derivative, it is possible to establish a fractional probabilistic Taylor's theorem in the Caputo sense too. The Caputo derivative, denoted by $_{*}D_{a+}^{\alpha } $, is defined by exchanging the operators $ I_{a+}^{m-\alpha } $ and $ D^{m} $ in the classical definition (\ref{RLderivative}).
Taking the paper of Odibat et al. \cite{Odibat2007} as a starting point, the following theorem, which is in some sense the equivalent of Theorem \ref{ThTaylor}, can be effortlessly proved.
\begin{theorem}
Let $ \alpha\in \left(0,1\right] $ and let $ X $ be a nonnegative random variable with distribution $ F $ and moment $ \mathbb{E}\left[X^{(n+1)\alpha}\right]<+\infty $ for some integer $ n\geq 0 $.
Assume that $g$ is a function defined on $I_{F}$, with $I_{F}$ defined in (\ref{eq:IF}), and suppose that $_{*}D_{0+}^{\alpha }g(x)\in C \left(I_{F}\right) $ for $ k=0,1,\dots, n+1 $ and $ \mathbb{E}\left [ \left | _{*}D_{0}^{\left ( n+1 \right )\alpha }g\left ( X_{\,\alpha }^{\left ( n+1 \right )} \right ) \right |\right ]<+\infty $. Then $ \mathbb{E}\left [ g\left ( X \right ) \right ]<+\infty $ and
\begin{multline*}
\mathbb{E}\left [ g(X) \right ]= \sum_{i=0}^{n}\frac{\left(_{*}D_{0+}^{i\alpha }f\right)(0)}{\Gamma \left ( i\alpha +1 \right )}\mathbb{E}\left [ X^{i\alpha } \right ]\\
+\frac{\mathbb{E}[X^{(n+1)\alpha}]}{\Gamma \left ( \left ( n+1 \right )\alpha +1 \right )}\,\mathbb{E}\left [ _{*}D_{0}^{\left ( n+1 \right )\alpha }g\left ( X_{\,\alpha }^{\left ( n+1 \right )} \right ) \right ],
\end{multline*}
where $ _{*}D_{0+}^{n\alpha }\,=\,_{*}D_{0+}^{\alpha }\cdot\,_{*}D_{0+}^{\alpha }\,\cdots\,_{*}D_{0+}^{\alpha } $ ($ n $ times) and $ X_{\,\alpha }^{\left ( n+1 \right )} $ has density (\ref{FractionalDensity}).
\end{theorem}

\section*{Acknowledgements}
The research that led to the present paper was partially supported by  a grant of the group GNCS of INdAM. 


\end{document}